# Some experimental observations about Hankel determinants of convolution powers of Catalan numbers

Johann Cigler

Fakultät für Mathematik Universität Wien

johann.cigler@univie.ac.at

**Abstract.**

Computer experiments suggest some conjectures about Hankel determinants of convolution powers of Catalan numbers. Unfortunately, for most of them I have no proofs. I would like to present them anyway hoping that someone finds them interesting and can prove them.

## 1. Some general observations

Let $C(x)=\frac{1-\sqrt{1-4x}}{2x}$ be the generating function of the Catalan numbers $C_n=\frac{1}{n+1}\binom{2n}{n}$.

The coefficients $C_{k,n}$ of $C(x)^k$ for a positive integer $k$ are the $k-$th convolution powers of the Catalan numbers. It is well known (cf. [4], (5.70)) that

$$C_{k,n}=\binom{2n+k}{n}\frac{k}{2n+k}=\frac{k}{n+k}\binom{2n+k-1}{n}=\binom{2n+k-1}{n}-\binom{2n+k-1}{n-1}. \tag{1}$$

If we extend $C_{k,n}$ to all $n\in\mathbb{Z}$ by setting $C_{k,n}=0$ for $n<0$ then the Hankel determinants

$$D_{k,m}(n)=\det\left(C_{k,i+j+m}\right)_{i,j=0}^{n-1} \tag{2}$$

are defined for $m\in\mathbb{Z}$ and positive integers $k,n.$ As usual we set $D_{k,m}(0)=1.$

Let

$$\mathbf{D}_{k,m}=\left(D_{k,m}(0),D_{k,m}(1),D_{k,m}(2),D_{k,m}(3),\cdots\right) \tag{3}$$

denote the sequence of the Hankel determinants of the shifted sequence $\left(C_{k,m+n}\right)_{n\geq 0}$ and let

$$T_r\mathbf{D}_{k,m}=\left(1,0,\cdots,0,D_{k,m}(0),D_{k,m}(1),D_{k,m}(2),\cdots\right) \tag{4}$$

be the sequence $\mathbf{D}_{k,m}$ translated $r$ units to the right with the first $r-1$ terms equal to $1,0,\cdots,0.$

---

As shown in [1] for $m \in \mathbb{N}$

$$\mathbf{D}_{1,-m} = (-1)^{\binom{m+1}{2}} T_{m+1} \mathbf{D}_{1,m+1}. \tag{5}$$

For $m = 0$ this reduces to the well-known fact

$$D_{1,0}(n) = D_{1,1}(n) = 1. \tag{6}$$

For the general case we obtained in [1]

**Conjecture 1**

*For $m \in \mathbb{N}$ we get*

$$\mathbf{D}_{2k,1-k-m} = (-1)^{\binom{m+k}{2}} T_{m+k} \mathbf{D}_{2k,1-k+m} \tag{7}$$

*and*

$$\mathbf{D}_{2k-1,2-k-m} = (-1)^{\binom{m+k-1}{2}} T_{m+k-1} \mathbf{D}_{2k-1,-k+m+1}. \tag{8}$$

For example,

$$\mathbf{D}_{3,1} = \left(1,3,3,-1,-6,-6,1,9,9,-1,-12,-12,\cdots\right),$$
$$\mathbf{D}_{3,-2} = -T_3 \mathbf{D}_{3,1} = \left(1,0,0,-1,-3,-3,1,6,6,-1,-9,-9,\cdots\right),$$

$$\mathbf{D}_{4,1} = \left(1,4,-4,-20,9,56,-16,-120,25,220,-36,-364,\cdots\right),$$
$$\mathbf{D}_{4,-3} = T_4 \mathbf{D}_{4,1} = \left(1,0,0,0,1,4,-4,-20,9,56,-16,-120,25,\cdots\right).$$

The special case $m = 0$ of (7) gives

$$\mathbf{D}_{2k,1-k} = (-1)^{\binom{k}{2}} T_k \mathbf{D}_{2k,1-k} \tag{9}$$

and therefore reduces to

**Theorem 2** ([2], Cor. 15)

$$\begin{aligned} &D_{2k,1-k}\left(kn\right) = (-1)^{n\binom{k}{2}}, \\ &D_{2k,1-k}(n) = 0 \quad \text{else}. \end{aligned} \tag{10}$$

The special case $m = 0$ of (8) gives

$$\begin{aligned} \mathbf{D}_{2k-1,2-k} &= (-1)^{\binom{k-1}{2}} T_{k-1} \mathbf{D}_{2k-1,1-k}, \\ \mathbf{D}_{2k-1,1-k} &= (-1)^{\binom{k}{2}} T_k \mathbf{D}_{2k-1,2-k}. \end{aligned} \tag{11}$$

By changing $k \to k+1$ this reduces to

**Theorem 3** ([2], Cor. 12 and Cor. 13)

*For* $k \geq 1$

$$\begin{aligned} &D_{2k+1,-k}\left((2k+1)n\right) = (-1)^{kn}, \\ &D_{2k+1,-k}\left((2k+1)n+k+1\right) = (-1)^{kn+\binom{k+1}{2}}, \\ &D_{2k+1,-k}\left((2k+1)n+j\right) = 0 \quad \text{else.} \end{aligned} \tag{12}$$

$$\begin{aligned} &D_{2k+1,1-k}\left((2k+1)n\right) = (-1)^{kn}, \\ &D_{2k+1,1-k}\left((2k+1)n+k\right) = (-1)^{kn+\binom{k}{2}}, \\ &D_{2k+1,1-k}\left((2k+1)n+j\right) = 0 \quad \text{else.} \end{aligned} \tag{13}$$

For example,

$$\begin{pmatrix} D_{3,-1}(n) \\ D_{3,0}(n) \end{pmatrix}_{0\leq n\leq 8} = \begin{pmatrix} 1 & 0 & -1 & -1 & 0 & 1 & 1 & 0 & -1 \\ 1 & 1 & 0 & -1 & -1 & 0 & 1 & 1 & 0 \end{pmatrix},$$

$$\begin{pmatrix} D_{5,-2}(n) \\ D_{5,-1}(n) \end{pmatrix}_{0\leq n\leq 9} = \begin{pmatrix} 1 & 0 & 0 & -1 & 0 & 1 & 0 & 0 & -1 & 0 \\ 1 & 0 & -1 & 0 & 0 & 1 & 0 & -1 & 0 & 0 \end{pmatrix}.$$

**Remark**

These facts suggest that from the point of view of Hankel determinants it would be more convenient to consider the shifted sequences $\left(c_{2k,n}\right)_{n\geq 0} = \left(C_{2k,n+1-k}\right)_{n\geq 0}$ and $\left(c_{2k-1,n}\right)_{n\geq 0} = \left(C_{2k-1,n+1-k}\right)_{n\geq 0}$ for $k \geq 1$.

## 2. Associated polynomials

It is well known (cf. e.g. [5], (5.67)) that for $m \in \mathbb{N}$ the Hankel determinants

$$D_{1,m}(n) = \det\left(C_{i+j+m}\right)_{i,j=0}^{n-1} \tag{14}$$

are polynomials $p_m(n)$

$$D_{1,m}(n) = p_m(n) = \prod_{1\leq i\leq j\leq m-1} \frac{2n+i+j}{i+j}. \tag{15}$$

An analog for $k > 1$ is

**Conjecture 4**

*For* $k \geq 1$, $j \in \{0,1,\cdots,k-1\}$ and $m \in \mathbb{Z}$

$$p_{2k,m,j}(n)=(-1)^{n\binom{k}{2}}D_{2k,m}(kn+j) \tag{16}$$

*and for* $j\in\{0,1,\cdots,2k\}$

$$p_{2k+1,m,j}(n)=(-1)^{nk}D_{2k+1,m}\left((2k+1)n+j\right) \tag{17}$$

*are polynomials in n.*

**2.1. The polynomials** $p_{2k,m,j}$.

**Conjecture 5**

*For* $n,k\geq 1$ *and* $m\geq -1$

$$\begin{aligned}&(-1)^{n\binom{k+m}{2}}D_{2k+2m,-m}((k+m)n)=(n+1)^{k-1},\\&(-1)^{n\binom{k+m}{2}}D_{2k+2m,-m}((k+m)n+1+m)=(-1)^{\binom{m+1}{2}}(n+1)^{k-1}.\end{aligned} \tag{18}$$

For example for $k=3$

$\left(D_{4,1}(n)\right)_{n\geq 0}=(\mathbf{1},-4,\mathbf{-4},-20,\mathbf{9},56,\mathbf{-16},-120,\mathbf{25},220,\mathbf{-36},-364,\mathbf{49},\cdots)$

$\left(D_{6,0}(n)\right)_{n\geq 0}=(\mathbf{1},1,-9,\mathbf{-4},-4,45,\mathbf{9},9,-126,\mathbf{-16},-16,270,\mathbf{25},\cdots)$,

$\left(D_{8,-1}(n)\right)_{n\geq 0}=(\mathbf{1},0,-1,-16,\mathbf{4},0,-4,-80,\mathbf{9},0,-9,-224,\mathbf{16},0,-16,\cdots)$.

$\left(D_{10,-2}(n)\right)_{n\geq 0}=(\mathbf{1},0,0,-1,25,\mathbf{4},0,0,-4,125,\mathbf{9},0,0,-9,350,\mathbf{16},0,0,-16,\cdots)$

**Conjecture 6** (Wang-Xin [6], Conjecture 10)

$$\begin{aligned}&\deg p_{2k,0,j}(n)=(2j-1)(k-j)\quad \text{for } k\geq 2j-1,\\&\textit{and}\quad \deg p_{2k,0,k+1-j}(n)=\deg p_{2k,0,j}(n).\end{aligned} \tag{19}$$

For example

| $2k$ | $j=2$ | $j=3$ | $j=4$ | $j=5$ | $j=6$ | $j=7$ | $j=8$ |
|---|---|---|---|---|---|---|---|
| 6 | 3 | | | | | | |
| 8 | 6 | 6 | | | | | |
| 10 | 9 | 10 | 9 | | | | |
| 12 | 12 | 15 | 15 | 12 | | | |
| 14 | 15 | 20 | 21 | 20 | 15 | | |
| 16 | 18 | 25 | 28 | 28 | 25 | 18 | |
| 18 | 21 | 30 | 35 | 36 | 35 | 30 | 21 |

For $k=3$ we get

$$p_{6,0,2}(n)=-3^2\frac{(n+1)(n+2)(2n+3)}{6}. \qquad (20)$$

Moreover for $m\geq -2$

$$p_{2(3+m),-m,2+m}(n)=(-1)^{\binom{m+2}{2}}(3+m)^2\frac{(n+1)(n+2)(2n+3)}{6}. \qquad (21)$$

For $k>3$ there is no simple relation between the polynomials $p_{2k,0,j}(n)$ and $p_{2(k+m),-m,j+m}(n)$, but there is one between their leading coefficients. Let us first collect some information about the leading coefficients of $p_{2k,0,j}(n)$.

For $p(n)=a_0+a_1n+\cdots+a_kn^k$ the leading coefficient $a_k$ is given by $\frac{\Delta^k}{k!}p(n)$ for each $n$ where $\Delta$ denotes the difference operator $\Delta p(n)=p(n+1)-p(n)$.

**Conjecture 7**

*For $j\geq 1$ the leading coefficient of $p_{2k,0,j}(n)$ is $A_{2k,j}k^{2(j-1)(k-j)}$ for some $A_{2k,j}$.*

Let us consider $A_{2k,j}$ for the first cases

| $2k$ | $j=2$ | $j=3$ | $j=4$ | $j=5$ | $j=6$ |
|---|---|---|---|---|---|
| 6 | $-\frac{1}{3}$ | | | | |
| 8 | $\frac{1}{45}$ | $-\frac{1}{45}$ | | | |
| 10 | $-\frac{2}{945}$ | $\frac{1}{4725}$ | $\frac{2}{945}$ | | |
| 12 | $\frac{1}{4725}$ | $\frac{1}{4465125}$ | $-\frac{1}{4465125}$ | $\frac{1}{4725}$ | |
| 14 | $-\frac{2}{93555}$ | $\frac{1}{2210236875}$ | $-\frac{1}{46414974375}$ | $-\frac{1}{2210236875}$ | $-\frac{2}{93555}$ |

With $\Phi_n=\prod_{j=1}^{n-1}(2j-1)!!$ the first terms $A_{2k,j}$ are integral multiples of $\frac{1}{\Phi_k}$ as shown in the following table.

$$\left(A_{2k,j}\right)_{3\le k\le 7, 2\le j\le 6} =$$

| $2k$ | $j=2$ | $j=3$ | $j=4$ | $j=5$ | $j=6$ |
|---|---|---|---|---|---|
| 6 | $-\frac{1}{\Phi_3}$ | | | | |
| 8 | $\frac{1}{\Phi_4}$ | $-\frac{1}{\Phi_4}$ | | | |
| 10 | $-\frac{10}{\Phi_5}$ | $\frac{1}{\Phi_5}$ | $\frac{10}{\Phi_5}$ | | |
| 12 | $-\frac{945}{\Phi_6}$ | $\frac{1}{\Phi_6}$ | $-\frac{1}{\Phi_6}$ | $\frac{945}{\Phi_6}$ | |
| 14 | $-\frac{992250}{\Phi_7}$ | $\frac{21}{\Phi_7}$ | $-\frac{1}{\Phi_7}$ | $-\frac{21}{\Phi_7}$ | $-\frac{992250}{\Phi_7}$ |

**Remark**

It is interesting to note that $\frac{1}{\Phi_m}$ is also the leading coefficient of $p_m(n)=\prod_{1\le i\le j\le m-1}\frac{2n+i+j}{i+j}$,

because $a(m):=\prod_{1\le i\le j\le m-1}(i+j)$ satisfies

$$a(m)=\prod_{1\le i\le j\le m-1}(i+j)=a(m-1)\prod_{1\le i\le m-1}(i+m-1)=a(m-1)\frac{(2m-2)!}{(m-1)!}=a(m-1)2^{m-1}(2m-1)!!.$$

For $j=2$ and $j=k-1$ we get by comparing $A_{2k,2}$ and $A_{2k,k-1}$ with

$$x\cot x=\sum_{k\ge 0}(-)^k\frac{2^{2k}B_{2k}}{(2k)!}x^{2k}=1-\frac{x^2}{3}-\frac{x^4}{45}-\frac{2x^6}{945}-\frac{x^8}{4725}-\frac{2x^{10}}{93555}-\frac{1382x^{12}}{638512875}-\cdots$$

**Conjecture 8**

*For* $k\ge 3$

$$(-1)^{n\binom{k}{2}}D_{2k,0}\left(kn+2\right) \tag{22}$$

*and*

$$(-1)^{n\binom{k}{2}+\binom{k+2}{2}}D_{2k,0}\left(kn+k-1\right) \tag{23}$$

*are polynomials of degree* $3(k-2)$ *with leading coefficient*

$$-\frac{B_{2k-4}}{(2k-4)!}(2k)^{2k-4}. \tag{24}$$

where $\left(B_{2n}\right)_{n\ge 0}=\left(1,\frac{1}{6},-\frac{1}{30},\frac{1}{42},-\frac{1}{30},\frac{5}{66},-\frac{691}{2730},\frac{7}{6},\cdots\right)$ are the Bernoulli numbers.

**Conjecture 9**

*Let* $k \geq 3,\ 2 \leq j \leq k-1$ *and* $m \geq -2.$ *If the leading coefficient of* $p_{2k,0,j}(n)$ *is* $A_{2k,j}k^{2(j-1)(k-j)}$ *then the leading coefficient of* $p_{2(k+m),-m,j+m}(n)$ *is* $(-1)^{\binom{m}{2}} A_{2k,j}(k+m)^{2(j-1)(k-j)}$ *for even* $j$ *and* $(-1)^{\binom{m+1}{2}} A_{2k,j}(k+m)^{2(j-1)(k-j)}$ *for odd* $j.$

**2.2. The polynomials** $p_{2k+1,m,j}(n)$.

**Conjecture 10**

$$\begin{aligned} &D_{2k+1,0}\left((2k+1)n\right) = D_{2k+1,0}\left((2k+1)n+1\right) = (-1)^{kn}, \\ &D_{2k+1,0}\left((2k+1)n+k+1\right) = 0, \end{aligned} \tag{25}$$

For example,

$$\left(D_{3,0}(n)\right)_{n\geq 0} = \left(1,1,0,-1,-1,0,1,1,0,-1,-1,0,1,1,0,-1,-1,0,1,1,0,-1,-1,\cdots\right),$$

$$\left(D_{5,0}(n)\right)_{n\geq 0} = \left(1,1,-5,0,5,1,1,-10,0,10,1,1,-15,0,15,1,1,-20,0,20,1,1,\cdots\right),$$

$$\left(D_{7,0}(n)\right)_{n\geq 0} = \left(1,1,-14,-49,0,49,329,-1,-1,-315,196,0,-196,-1687,1,1,\cdots\right).$$

**Conjecture 11**

*For* $0 \leq m \leq k+1$

$$D_{2k+1,m-k+1}\left((2k+1)n+k\right) = (-1)^{kn+\binom{k}{2}}(2k+1)^m(n+1)^m.$$

For example,

$$\left(D_{3,0}(3n+1)\right)_{n\geq 0} = \left(1,-1,1,-1,1,-1,1,-1,1,\cdots\right),$$

$$\left(D_{3,1}(3n+1)\right)_{n\geq 0} = \left(3,-6,9,-12,15,-18,21,-24,27,\cdots\right),$$

$$\left(D_{3,2}(3n+1)\right)_{n\geq 0} = \left(9,-36,81,-144,225,-324,441,-576,729,\cdots\right).$$

**Conjecture 12** (Wang-Xin[6], Conjecture 9)

*For* $k \geq 1$ *and* $2 \leq j \leq k-1$

$(-1)^{nk} D_{2k+1,0}\left((2k+1)n+j\right)$ *and* $(-1)^{nk} D_{2k+1,0}\left((2k+1)n+2k+2-j\right)$

*are polynomials in* $n$ *of degree* $(j-1)\left(2k+1-2j\right).$

The following table shows the degrees. The degree of the polynomial $0$ is chosen to be $-1.$

| $2k+1$ | $j=2$ | $j=3$ | $j=4$ | $j=5$ | $j=6$ | $j=7$ | $j=8$ | $j=9$ | $j=10$ |
|---|---|---|---|---|---|---|---|---|---|
| 3 | $-1$ | | | | | | | | |
| 5 | 1 | $-1$ | 1 | | | | | | |
| 7 | 3 | 2 | $-1$ | 2 | 3 | | | | |
| 9 | 5 | 6 | 3 | $-1$ | 3 | 6 | 5 | | |
| 11 | 7 | 10 | 9 | 4 | $-1$ | 4 | 9 | 10 | 7 |

**Conjecture 13**

*For* $2 \le j \le k$ *the leading term of* $p_{2k+1,0,j}(n)$ *is* $B_{2k+1,j}(2k+1)^{(j-1)(1+2(k-j))}$ *for some numbers* $B_{2k+1,j}$. *For* $k < j \le 2k$ $B_{2k+1,0,j} = (-1)^{j-k-1} B_{2k+1,0,2k+2-j}$.

Table $\left(B_{2k+1,j}\right)_{k\ge1,j\ge2}$:

| $2k+1$ | $j=2$ | $j=3$ | $j=4$ | $j=5$ | $j=6$ | $j=7$ | $j=8$ | $j=9$ | $j=10$ | $j=11$ | $j=12$ |
|---|---|---|---|---|---|---|---|---|---|---|---|
| 3 | 0 | | | | | | | | | | |
| 5 | $-1$ | 0 | 1 | | | | | | | | |
| 7 | $\frac{1}{3}$ | $-1$ | 0 | 1 | $\frac{1}{3}$ | | | | | | |
| 9 | $-\frac{2}{15}$ | $\frac{1}{45}$ | 1 | 0 | $-1$ | $\frac{1}{45}$ | $\frac{2}{15}$ | | | | |
| 11 | $\frac{17}{315}$ | $\frac{1}{4725}$ | $-\frac{2}{945}$ | 1 | 0 | $-1$ | $-\frac{2}{945}$ | $-\frac{1}{4725}$ | $\frac{17}{315}$ | | |
| 13 | $-\frac{62}{2835}$ | $\frac{1}{297675}$ | $-\frac{1}{4465125}$ | $-\frac{1}{4725}$ | $-1$ | 0 | 1 | $-\frac{1}{4725}$ | $\frac{1}{4465125}$ | $\frac{1}{297675}$ | $\frac{62}{2835}$ |

Comparing with

$$\tan(x) = x + \frac{1}{3}x^3 + \frac{2}{15}x^5 + \frac{17}{315}x^7 + \frac{62}{2835}x^9 + \frac{1382}{155925}x^{11} + \cdots = \sum_{k\ge1} (-1)^{k-1} \frac{2^{2k}\left(2^{2k}-1\right)B_{2k}}{(2k)!} x^{2k-1}$$

and with $x\cot x$ we get

**Conjecture 14**

*The leading coefficient of* $(-1)^{k-1} p_{2k+1,0,2}(n)$ *and of* $p_{2k+1,0,2k}$ *is*

$$\frac{2^{2k-2}\left(2^{2k-2}-1\right)B_{2k-2}}{(2k-2)!}(2k+1)^{2k-3}. \qquad (26)$$

*The leading coefficient of* $p_{2k+1,0,k-1}(n)$ *is*

$$(-1)^{\binom{k-2}{2}}\frac{2^{2k-4}B_{2k-4}}{(2k-4)!}(2k+1)^{3k-6}. \tag{27}$$

**3. Generating functions**

Let us also give some conjectures for the generating functions of the Hankel determinants $D_{k,m}(n)$.

Let us first recall some definitions. A polynomial $P(x)$ of degree $d$ is called palindromic if $P\left(\frac{1}{x}\right)x^d = P(x)$ and skew palindromic if $P\left(\frac{1}{x}\right)x^d = -P(x)$.

Let $\pi_k(n)$, $k \geq 3$, denote the so called second k-gonal numbers $\pi_k(n) = \frac{n}{2}\big((k-2)n + (k-4)\big)$.

Note that $\pi_3(n) = \binom{n}{2}$, $\pi_4(n) = n^2$, $\pi_5(n) = \frac{n(3n+1)}{2}$.

**Conjecture 15**

*For* $k \geq 1$ *and* $m \geq 1-k$

$$\left(1-(-1)^{\binom{k}{2}}x^k\right)^{\binom{k+m}{2}+1}\sum_{n\geq 0} D_{2k,m}(n)x^n = P_{2k,m}(x) \tag{28}$$

*with* $\deg P_{2k,m}(x) = \pi_{k+2}(m+k-1)$.

*All polynomials* $P_{2k,m}(x)$ *for* $m \geq 1-k$ *are either palindromic or skew palindromic.*

*If* $k \equiv 1 \bmod 4$ $P_{2k,m}(x)$ *is palindromic for all* $m$, *if* $k \equiv 0 \bmod 4$ *then it is palindromic for odd* $m$ *and skew palindromic for even* $m$.

*For* $k \equiv 2 \bmod 4$ *it is palindromic for* $m \equiv 0 \bmod 4$ *and* $m \equiv 3 \bmod 4$, *for* $k \equiv 3 \bmod 4$ *it is palindromic if* $m \equiv 1 \bmod 4$ *and* $m \equiv 2 \bmod 4$.

For $k = 1$ this reduces to

$$(1-x)^{\binom{m+1}{2}+1}\sum_{n\geq 0} D_{2,m}(n)x^n = (1-x)^{\binom{m+1}{2}+1}\sum_{n\geq 0} D_{1,m+1}(n)x^n = P_{2,m}(x) \tag{29}$$

with $\deg P_{2,m}(x) = \binom{m}{2}$.

For example, we get $P_{2,0}(x)=P_{2,1}(x)=1,\ P_{2,2}(x)=1+x,\ P_{2,3}(x)=1+7x+7x^2+x^3,$
$P_{2,4}(x)=1+31x+187x^2+330x^3+187x^4+31x^5+x^6.$

For $k=2$ we get $\deg P_{4,m}(x)=(m+1)^2$. The first corresponding polynomials are

$P_{4,-1}(x)=1,\ P_{4,0}(x)=1+x,\ P_{4,1}(x)=1+4x-4x^3-x^4,$
$P_{4,2}(x)=1+14x+13x^2-111x^3-119x^4+119x^5+111x^6-13x^7-14x^8-x^9,$

$P_{4,3}(x)=1+48x+242x^2-1760x^3-7960x^4+10112x^5+47918x^6-9680x^7-84370x^8$
$-9680x^9+47918x^{10}+10112x^{11}-7960x^{12}-1760x^{13}+242x^{14}+48x^{15}+x^{16}.$

For $k=3$ we get

$P_{6,-2}(x)=1,\ P_{6,-1}(x)=1-x^2,\ P_{6,0}(x)=1+x-9x^2+9x^5-x^6-x^7,$

$P_{6,1}(x)=1+6x-69x^2-x^3+63x^4+561x^5-8x^6-609x^7-609x^8-8x^9+561x^{10}+63x^{11}-x^{12}-69x^{13}+6x^{14}+x^{15}.$

More generally $P_{2k,1-k}(x)=1$ for $k\geq 1$, $P_{2k,2-k}(x)=1+(-1)^{\binom{k-1}{2}}x^{k-1}$ for $k\geq 2$ and

$$P_{2k,3-k}(x)=\left(1+(-1)^{\binom{k-2}{2}}x^{k-2}\right)\left(1-x^{2k}\right)+(-1)^{\binom{k-1}{2}}\left(1+(-1)^{\binom{k}{2}}x^{k}\right)k^2x^{k-1}\ \text{ for } k\geq 3,$$

$$\left(P_{2k,3-k}(x)\right)_{k\geq 3}=\left(1+x-9x^2+9x^5-x^6-x^7,1-x^2-16x^3-16x^7-x^8+x^{10},1-x^3+25x^4+25x^9-x^{10}+x^{13},\cdots\right).$$

**Conjecture 16**

*For* $k\geq 1$ *and* $m\geq 1-k$

$$\left(1+(-1)^k x^{2k-1}\right)^{\binom{k+m}{2}+1}\sum_{n\geq 0}D_{2k-1,m}(n)x^n=Q_{2k-1,m}(x) \tag{30}$$

*with* $\deg Q_{2k-1,m}(x)=\pi_{2k+1}\left(m+k-1\right)+k.$

*All polynomials* $Q_{2k-1,m}(x)$ *are either palindromic or skew palindromic.*

*For* $k\equiv 1 \bmod 4$ $Q_{2k-1,m}(x)$ *is palindromic for even* $m$, *for* $k\equiv 3 \bmod 4$ *for odd* $m$,

*For* $k\equiv 0,2 \bmod 4$ $Q_{2k-1,m}(x)$ *is palindromic for* $m\equiv 0,1 \bmod 4.$

For $k=1$ (30) is not optimal. Comparing (29) with (30) gives $Q_{1,m+1}(x)=(x-1)^{m+1}P_{2,m}(x).$

For $k=2$ $\deg Q_{3,m}(x)=\pi_5(m+1)+2$. The first terms are $2,4,9,17,28,\cdots$.

$Q_{3,-1}(x)=1-x^2,\ Q_{3,0}(x)=(1+x)^2(1-x+x^2),\ Q_{3,1}(x)=(1+x)^5(1-x+x^2)^2,$

$Q_{3,2}(x)=(1-x)(1+x)^8(1-x+x^2)^3(1+5x+x^2).$

More generally

$$Q_{2k-1,1-k}(x)=1+(-1)^{\binom{k}{2}}x^k \text{ for } k\ge 2,$$

$$Q_{2k-1,2-k}(x)=\left(1+(-1)^{\binom{k-1}{2}}x^{k-1}\right)\left(1+(-1)^k x^{2k-1}\right) \textit{ for } k\ge 2.$$

$$Q_{2k-1,3-k}(x)=\left(1+(-1)^{\binom{k-2}{2}}x^{k-2}\right)\left(1+(-1)^k x^{2k-1}\right)+(2k-1)x^{k-1}\left(x^{k-1}+(-1)^{\binom{k-1}{2}}\right) \text{ for } k\ge 3.$$

## 4. Similar results for Hankel determinants of certain binomial coefficients

Similar facts hold for the sequences $\left(c_{k,n}\right)_{n\ge 0}$ with $c_{k,n}=\binom{2n+k}{n}$ for $n\ge 0$ and $c_{k,n}=0$ for $n<0$. Let $d_{k,m}(n)=\det\left(c_{k,i+j+m}\right)_{i.j=0}^{n-1}$.

Some results of [3] for $m=0$ suggest

**Conjecture 17.**

*For* $m\ge -1$

$$d_{2k+2m+1,-m}\left((2k+2m+1)n\right)=(2n+1)^k,$$

$$d_{2k+2m+1,-m}\left((2k+2m+1)n+1+m\right)=(-1)^{\binom{m+1}{2}}(2n+1)^k,$$

$$d_{2k+2m+1,-m}\left((2k+2m+1)n+k+m+1\right)=(-1)^{\binom{m+k+1}{2}}4^k(n+1)^k,$$

$$d_{2k+2m,-m}\left((2k+2m)n\right)=(-1)^{(k+m)n},$$

$$d_{2k+2m,-m}\left((2k+2m)n+m+1\right)=(-1)^{\binom{m+1}{2}+(k+m)n}.$$

An interesting analog of Conjectures 8 and 14 is

**Conjecture 18**

*For* $m \ge 0$ *the leading coefficients* $a_{2k+2m,-m,2k-3}$ *of the polynomials* $q_{2k+2m,-m}(n) = d_{2k+2m,-m}((2k+2m)n+m+2)(-1)^{(k+m)n}$ *of degree* $2k-3$ *are*

$$a_{2k+2m,-m,2k-3} = (-1)^{\binom{m}{2}+1} \frac{(k+m)^{2k-3} 4^{2k-2} \left(2^{2k-2}-1\right) B_{2k-2}}{(2k-2)!}$$

*and the leading coefficients* $a_{2k+2m+1,-m,3k-3}$ *of the polynomials* $q_{2k+2m+1,-m,3k-3}(n) = d_{2k+2m+1,-m}((2k+2m+1)n+2+m)$ *of degree* $3k-3$ *are*

$$a_{2k+2m+1,-m,3k-3} = (-1)^{\binom{m}{2}+1} \frac{2^{3k-2}(2k+m+1)^{2k-2} B_{2k-2}}{(2k-2)!}.$$